\documentclass[reqno,centertags, 12pt]{amsart}

\usepackage{amsmath,amsthm,amscd,amssymb,mathtools}
\usepackage{esint}  
\usepackage{latexsym}
\usepackage{graphicx}
\usepackage{enumitem}
\usepackage[mathscr]{eucal}

\usepackage{hyperref}
\newcommand{\bbC}{{\mathbb{C}}}

\newcommand{\bbE}{{\mathbb{E}}}

\newcommand{\bbR}{{\mathbb{R}}}

\newcommand{\bbZ}{{\mathbb{Z}}}

\newcommand{\calG}{{\mathcal G}}
\newcommand{\calH}{{\mathcal H}}

\newcommand{\calT}{{\mathcal T}}

\newcommand{\lb}{\label}

\newcommand{\spec}{\text{\rm{spec}}}

\newcommand{\bi}{\bibitem}

\newcommand{\beq}{\begin{equation}}
\newcommand{\eeq}{\end{equation}}
\newcommand{\ba}{\begin{align}}
\newcommand{\ea}{\end{align}}

\newcounter{smalllist}

\newcommand{\comm}[1]{}

\DeclareMathOperator{\Ima}{Im}

\allowdisplaybreaks
\numberwithin{equation}{section}

\newtheorem{theorem}{theorem}[section]

\theoremstyle{definition}
\newtheorem*{remark}{remark}

\newcommand{\jap}[1]{\langle #1 \rangle}

\newcounter{pnas}
\begin{document}

\title[Useful Formula]{A Useful Formula for Periodic Jacobi Matrices on Trees}
\author[J.~Banks, J.~Breuer, J.~Garza~Vargas, E.~Seelig and B.~Simon]{Jess Banks$^{a,1}$, Jonathan Breuer$^{b,2,3}$,\\Jorge Garza Vargas$^{c,2,4}$, Eyal Seelig$^{d,2,3}$ \\and Barry Simon$^{e,2}$}

\thanks{$^a$ Department of Mathematics, University of California–Berkeley, Berkeley, USA. Email: jess.m.banks@berkeley.edu}

\thanks{$^b$ Institute of Mathematics, The Hebrew University of Jerusalem, Jerusalem, 91904, Israel. E-mail: jbreuer@math.huji.ac.il.}

\thanks{$^c$ Department of Computing and Mathematical Sciences, MC 305-16, California Institute of Technology, Pasadena, CA 91125, USA. E-mail: jgarzav@caltech.edu.}

\thanks{$^d$ Institute of Mathematics, The Hebrew University of Jerusalem, Jerusalem, 91904, Israel. E-mail: eyal.seelig@mail.huji.ac.il.}

\thanks{$^e$ Departments of Mathematics and Physics, Mathematics 253-37, California Institute of Technology, Pasadena, CA 91125, USA. E-mail: bsimon@caltech.edu.}

\thanks{$^1$ TK.}

\thanks{$^2$ Research of JBr, JGV, EY and BS supported in part by Israeli BSF Grant No. 2020027. }

\thanks{$^3$ Research of JBr and EY supported in part by Israel Science Foundation Grant No. 1378/20}

\thanks{$^4$ Research supported in part by NSF FRG Award 1952777 and Caltech Carver Mead New Adventures Fund under the aegis of Joel Tropp's award, and Caltech CMI Postdoctoral Fellowship }

\

\date{\today}

\begin{abstract}  We introduce a function of the density of states for periodic Jacobi matrices on trees and prove a useful formula for it.  This allows new, streamlined proofs of the gap labelling and Aomoto index theorems.  We prove a version of this new formula for the Anderson model on trees.
\end{abstract}

\maketitle

\section{Introduction} \lb{s1}

Graph Jacobi matrices provide a unified framework for dealing with graph adjacency matrices, weighted Laplacians and Schr\"odinger operators. Their spectral theory therefore has connections with various fields, among those are mathematical physics, analysis, probability and number theory. This note deals with periodic Jacobi matrices on trees, which arise through viewing the tree as the universal cover of a finite graph. Such operators have attracted a considerable amount of interest recently~\cite{AFH, AomotoPoint, ABKS, ABS, BGM, Cartier, CSZ, FTP, FTS, GVK, KLW1, KLW2, NW, Sun, WoessBk}. The purpose of this note is to announce and give an interim report on the use of a new formula which, in particular, provides a short proof of Sunada's gap labelling result~\cite{Sun}, without the use of $C^*$ algebras.

We start out with a connected, finite graph, $\calG$, which can have self-loops and multiple edges between a pair of vertices but which, for simplicity of exposition, we suppose is leafless.  We use $V(\calG)$ for the vertex set of $\calG$ and $E(\calG)$ (sometimes just $V$ and $E$) for the set of edges.  We pick an orientation for each edge, $e$, using $\check{e}$ for the oppositely directed edge. $\sigma(e)$ is the initial vertex and $\tau(e)$ the final of the directed edge $e$, so for example, $\sigma(\check{e})=\tau(e)$.
We let $\tilde{E}$ denote the set of all edges with arbitrary assigned orientation so that $\#(\tilde{E})=2\#(E)$. We assign a potential, $b(v)\in\bbR$, to each vertex and coupling, $a(e)=a(\check{e})>0$, to each edge, calling these the \emph{Jacobi parameters} of $\calG$.

Let $\calT$ be the universal cover of $\calG$ -  it is always an infinite tree, and let $\pi:\calT\to\calG$ be the covering map.  We can lift the Jacobi parameters of $\calG$ to $\calT$ by setting $b(\tilde{v}) = b(\pi(\tilde{v}));\quad a(\tilde{e}) = a(\pi(\tilde{e}))$.  One defines an infinite matrix, $H$, indexed by $V(\calT)$ by
\begin{equation}\label{1.1}
  H_{\tilde{v}\tilde{w}} = \left\{
                             \begin{array}{ll}
                               b(\tilde{v}), & \hbox{ if }\tilde{v}=\tilde{w} \\
                               a(\tilde{e}), & \hbox{ if } (\tilde{v}\tilde{w}) = \tilde{e} \text{ an edge in }\tilde{E}(\calT)\\
                               0, & \hbox{ otherwise}
                             \end{array}
                           \right.
\end{equation}
and a corresponding bounded self-adjoint operator, $H$, on $\calH=\ell^2(V(\calT))$.  One defines the \emph{period}, $p$, to be $\#(V(\calG))$.  If $\calG$ is a single cycle, then $\calT$ is $\bbZ$ and the Jacobi parameters are periodic in the naive sense.  This classical subject (of $1D$ periodic Jacobi matrices) has been extensively studied; see for example, Simon~\cite[Chaps. 5, 6, 8]{SiSz}.

Deck transformations induce unitary maps on $\calH$ which commute with $H$.  In particular, for every $v\in V(\calG)$, the spectral measure, $d\mu_{\tilde{v}}$, and Green's function, $\jap{\delta_{\tilde{v}},(H-z)^{-1}\delta_{\tilde{v}}}$, are the same for all $\tilde{v}\in V(\calT)$ with $\pi(\tilde{v})=v$.  We use $d\mu_v$ and $G_v(z)$ for these common values.  It is a basic fact that in one form goes back at least to \cite{Chom} (see also \cite{Woess, KLW1}) that each $G_v(z)$ defined for $z\in\bbC_+$ is an algebraic function which can be continued across the real axis with finitely many points removed (this implies, see \cite[Theorem 6.7]{ABS}, that the spectrum of $H$ has no singular continuous part and the densities of the a.c. part of the spectral measures are real analytic in the interior of the spectrum except for possible algebraic singularities).

One defines the \emph{density of states} measure, $dk(E)$ (and \emph{integrated density of states}, aka IDS, $k(E) = dk((-\infty,E))$), by
\begin{equation}\label{1.2}
  dk = \frac{1}{p} \sum_{v\in V} d\mu_v
\end{equation}
\begin{remark} For Jacobi matrices on $\bbZ^\nu$, the analog is the limit of number of eigenvalues in a large box (with, say, free boundary conditions) per unit volume;
	because truncated trees have so many boundary points, that is not true for trees with general boundary conditions (BC) although one can carefully choose periodic BC so that it is~\cite{ABKS}.
\end{remark}
The support of the measure $dk$ is the spectrum of $H$ and by the definition of spectral measures, one has that
\begin{equation}\label{1.3}
  \int \frac{1}{\lambda-z}\,dk(\lambda) = \frac{1}{p} \sum_{v\in V}G_v(z)
\end{equation}

One of the fundamental results in the theory is

\begin{theorem} [Sunada~\cite{Sun}] \label{T1.1} In any gap of the spectrum of $H$, the IDS is an integral multiple of $1/p$.  In particular, the spectrum has at most $p$ connected components.
\end{theorem}

Sunada's proof, while elegant, is involved since it uses some deep results of Pimsner-Voiculescu~\cite{PV} from the $K$-theory of $C^*$-algebras.  One of our main new results is a short proof of Sunada's theorem that, in particular, makes no use of $C^*$-algebras.

Another fundamental result is the Aomoto index theorem.  In the $1D$ case, $H$ does not have any point spectrum but in other cases that is not true - see, for example, Avni et al.~\cite[Example 7.2]{ABS} or the extensive study in Banks et al.~\cite{BGM}.  In that case, given an eigenvalue, $\lambda$, define $X_1(\lambda)$ to the set of vertices, $v\in V$, so that some $\tilde{v}$  with $\pi(\tilde{v})=v$ has some energy $\lambda$ eigenfunction, $\psi$, with $\psi(\tilde{v})\ne 0$. Define $\partial X_1(\lambda)$ to be those $v\in V$ not in $X_1(\lambda)$ but neighbors of points in $X_1(\lambda)$, and we let $E(\lambda)$ be the set of edges with both endpoints in $X_1(\lambda)$.

\begin{theorem} [Aomoto Index Theorem~\cite{AomotoPoint}] \label{T1.2} The measure $dk$ has a mass at an eigenvalue, $\lambda$, of weight $I(\lambda)/p$ where
\begin{equation}\label{1.4}
  I(\lambda) = \#(X_{1}(\lambda)) - \#(\partial X_1(\lambda)) - \#(E(\lambda))
\end{equation}
\end{theorem}

A second proof of this theorem can be found in Banks et al.~\cite{BGM}.  Both earlier proofs involve detailed combinatorial analyses.  The second of our new results here is a different proof of the Aomoto index theorem that some may find simpler but that, in any event, is very illuminating.

Our new approach concerns a basic function which we will call the \emph{Floquet function} defined in $\bbC_+$ by
\begin{equation}\label{1.5}
  \Phi(z) = \exp\left(p\int \log(t-z) \, dk(t)\right)
\end{equation}
which clearly has an analytic continuation to a neighborhood of $\bbC_+\cup(\bbR\setminus\spec(H))$.  In the $1D$ case, the Thouless formula (\cite[Theorem 5.4.12]{SiSz}) implies (if $\prod_{j=1}^{p}a_j=1$) that if $u_j(z)$ is a \emph{Floquet solution} (i.e. solution of the difference equation
\begin{equation}\label{1.5A}
  a_j u_{j+1}+b_j u_j+a_{j-1}u_{j-1} = zu_j
\end{equation}
with $u_{j+p} = A u_j$ for a constant $A$), then (\cite[Theorem 5.4.15]{SiSz}) $A=\Phi(z)$ or $\Phi(z)^{-1}$  which is why we give $\Phi$ this name.

Because of \eqref{1.3} we have that
\begin{equation}\label{1.6}
  \frac{d}{dz} \log(\Phi(z)) = -\sum_{v\in V} G_v(z)
\end{equation}

In Section \ref{s2}, we'll prove an explicit formula for the Floquet function in terms of Green's functions and $m$-functions (objects whose definition we recall there).  In Section \ref{s3}, we'll use this Floquet formula to prove the Sunada gap labelling theorem and in Section \ref{s4}, we'll sketch our new proof of the Aomoto index theorem (in the case where the eigenvalue is isolated from the continuous spectrum; see the discussion there).   In Section \ref{s5}, we will discuss a version of the Floquet formula for the Anderson model on trees.  Since, as we'll explain, one can regard the Floquet formula as half a Thouless formula, we hope to find some interesting applications of that result.

\section{The Floquet Formula} \lb{s2}

We will prove a useful formula for the Floquet function.  To do so, we need to recall what the $m$-functions are and the relations between the Green's and $m$-functions.  Given $e\in E$, pick $\tilde{e}\in E(\calT)$ with $\pi(\tilde{e})=e$.  Removing $\tilde{e}$ from $\calT$ breaks that graph into two pieces, $\calT^+_{\tilde{e}}$ with $\tau(\tilde{e})$ and $\calT^-_{\tilde{e}}$ with $\sigma(\tilde{e})$.  We let $H_{\tilde{e}}^\pm$ be the operators on $\ell^2(V(\calT^\pm_{\tilde{e}}))$ with the restricted Jacobi parameters and set
\begin{equation}\label{2.1}
  m_e(z) = \jap{\delta_{\tau(\tilde{e})},(H_{\tilde{e}}^+-z)^{-1}\delta_{\tau(\tilde{e})}}
\end{equation}
The use of deck transformations shows this depends only on $e$ and not the choice of $\tilde{e}$ over $e$.

The use of the method of Schur complements (see \cite[Section 6]{ABS} for a proof; the formulae appear at least as early as \cite[Proposition 2.1]{Klein}) shows that
\begin{align}
  \frac{1}{G_u(z)} &= -z+b_u-\sum_{\mathclap{f\in\tilde{E}:\;\sigma(f)=u}} a_f^2m_f(z) \label{2.2} \\
  \frac{1}{m_f(z)} &= -z+b_u-\sum_{\mathclap{\substack{f'\in\tilde{E},f'\ne \check{f}\\ \sigma(f')=\tau(f)}}} a_{f'}^2m_{f'}(z) \label{2.3}
\end{align}
which implies for any $e\in\tilde{E}$ that
\begin{equation}\label{2.4}
  G_{\sigma(e)} = \frac{1}{m_{\check{e}}^{-1}-a_e^2m_e} = \frac{m_{\check{e}}}{1-a_e^2m_em_{\check{e}}}
\end{equation}

Define
\begin{equation}\label{2.5}
  Q_e(z) = \frac{1}{1-a_e^2m_e(z)m_{\check{e}}(z)} = \frac{G_{\sigma(e)}(z)}{m_{\check{e}}(z)} = \frac{G_{\tau(e)}(z)}{m_{e}(z)}
\end{equation}

We are heading towards the proof of a lovely formula we call the \emph{Floquet formula}:

\begin{theorem} [Floquet Formula] \label{T2.1} We have that
\begin{equation}\label{2.6}
  \Phi(z) = \frac{\prod_{e\in E(\calG)}Q_e(z)}{\prod_{u\in V(\calG)}G_u(z)}
\end{equation}
initially for $z\in\bbC_+$, but the right side defines a meromorphic continuation to $(\bbC\setminus\spec(H))\cup(\text{isolated point spectrum of }H )$.
\end{theorem}

\begin{remark} Using \eqref{2.5}, this can also be written
\begin{equation}\label{2.6A}
  \Phi(z) = \frac{\prod_{e\in E(\calG)}G_{\tau(e)}(z)}{\prod_{u\in V(\calG)}G_u(z)\prod_{e\in E(\calG)}m_e(z)}
\end{equation}
In particular, this implies that $\Phi$ is an algebraic function.
\end{remark}

We sketch our proof of this result.  Let $\Psi$ be the right side of \eqref{2.6}.  It is easy to see that as $x\to\infty$ in $\bbR$, that $\Phi(-x) = x^p+\text{O}(x^{p-1})$ and $\Psi(-x) = x^p+\text{O}(x^{p-1})$ so to prove \eqref{2.6}, it suffices to prove that for $z\in\bbC_+$
\begin{equation}\label{2.7}
  \log(\Psi)'(z) = \log(\Phi)'(z)
\end{equation}
where $\cdot'=\tfrac{d\cdot}{dz}$.

To compute $ \log(\Psi)'(z)$, we note that, by \eqref{2.2}, we have that
\begin{align}\label{2.8}
  (\log(G_u))' &= -\left(\log\left(\frac{1}{G_u}\right)\right)'= -G_u\left(\frac{1}{G_u}\right)' \nonumber\\
               &= G_u+\sum_{\mathclap{e\in\tilde{E}:\;\sigma(e)=u}}a_e^2 m_e'G_u
\end{align}
and that by, \eqref{2.5},
\begin{align}\label{2.9}
  (\log(Q_e))' &= (a_e^2m_e'm_{\check{e}}+a_e^2m_em_{\check{e}}')Q_e \nonumber \\
               &= a_e^2G_{\sigma(e)}m_e'+a_e^2G_{\tau(e)}m_{\check{e}}'
\end{align}

Therefore
\begin{align}\label{2.10}
  \sum_{e\in E} (\log(Q_e))' &= \sum_{e\in\tilde{E}} a_e^2G_{\sigma(e)}m_e' \nonumber \\
                               &= \sum_{u\in V}\qquad\sum_{\mathclap{e\in\tilde{E}:\;\sigma(e)=u}}a_e^2 m_e'G_u \nonumber\\
                            &= \sum_{u\in V} (-G_u+(\log(G_u))')
\end{align}
which, by \eqref{1.6}, proves \eqref{2.7} and so Theorem \ref{T2.1}.

\section{Gap Labelling} \lb{s3}

In this section, we present our new proof of Sunada's Gap Labelling theorem, Theorem \ref{T1.1}.  Basically, it is an immediate consequence of the Floquet formula \eqref{2.6}.
We need some care in determining the branch of $\log$ used $\eqref{1.5}$.  We pick the branch where when $z\in\bbC_+$ is taken near $-\infty$ on the real axis, $\Phi$ has an argument near $0$.  That is, we are using the branch where when $z=-x$ ($x$ near $+\infty$) and $t$ in a bounded interval, we have that $\log(t-z)>0$ and we are then continuing $z$ through the upper plane.   Thus, if $E_0$ is a real point in the resolvent set of $H$, the integral defining $\Phi$, \eqref{1.5}, can be analytically continued from $\bbC_+$ to a neighborhood of $E_0$ and for $s=t-E_0\ne 0$ real, we have that
\begin{equation}\label{3.1}
  \Ima(\log(s)) = \left\{
                    \begin{array}{ll}
                      0, & \hbox{if } s>0 \\
                      -\pi, & \hbox{ if } s<0
                    \end{array}
                  \right.
\end{equation}
Moreover, the Floquet formula can be analytically continued to a set including $E_0$. Thus
\begin{equation}\label{3.2}
  \Ima\left(p\int\log(t-E_0)\,dk(t)\right) = - p\pi k(E_0)
\end{equation}

That means that $pk(E_0)\in\bbZ\iff \Phi(E_0)$ is real.  But for $x\in\bbR\setminus\spec(H)$, each $G_v(x)$ and $m_e(x)$ is analytic (meromorphic for $m$), we see that except for potential isolated poles (actually, it is easy to see that $\Phi$ has no poles), $\Phi$ is real in gaps!

\section{Aomoto Index Theorem} \lb{s4}

In this section, we will sketch (with full details in a later publication) a proof of the Aomoto Index Theorem (Theorem \ref{T1.2}) at least in the case where the eigenvalue is an isolated point of the spectrum (we hope in the later publication to deal with the general case; we'll explain the potential difficulty soon - see point (1) below; the next paragraph also uses that the eigenvalue is isolated).  We note that the earlier proofs of this theorem (\cite{AomotoPoint, BGM}) handle the general case and that Banks et al.~\cite{BGM} provide examples where there are non-isolated eigenvalues and also where there are isolated eigenvalues.

The Floquet function is involved with the question of the weight of an eigenvalue because, by the discussion in the last section, $\lambda$ is an isolated eigenvalue with $dk$-weight $I/p$ if and only if the argument of $\Phi(x)$ jumps by $I\pi$ as $x$ passes through $\lambda$.  For isolated eigenvalues, by the Sunada theorem, $I$ is an integer so this happens if and only if $\Phi$ has a zero of order $I$ at $\lambda$.

The punch line is that \eqref{1.4} will come from the Floquet formula, \eqref{2.6}, and the fact that $G_v(z)$ has a simple pole at $v$ if and only if $v\in X_1(\lambda)$, it has simple zeros when $v\in\partial X_1(\lambda)$ and $Q_e$ has a simple pole if and only if $e\in E$.  There can be some additional zeros of $G_v$ and $Q_e$ but we will see that they cancel.

We will use $X_0(\lambda)= V\setminus (X_1(\lambda)\cup\partial X_1(\lambda))$. Henceforth, without loss, we can suppose that $\lambda=0$ for simplicity of notation and we drop $(\lambda)$ from $X_{0,1}(\lambda)$.

The proof relies on a sequence of observations:

(1) If $0$ is an isolated point in the spectrum then all Green's and $m$-functions are meromorphic in a neighborhood of $0$. If they have poles they are simple with negative residue and if they are zero, the zeros are simple with positive derivative (this follows from the fact that by the spectral theorem, the derivative of Green's and $m$-functions away from poles are strictly positive).  Thus in counting the order of a zero in \eqref{2.6A}, each $G$ or $m$ contributes either a single $+1$, $-1$ or $0$.  (For non-isolated zeros, the functions are only algebraic and so have Laurent-Puiseux series -- one needs to track potential fractional powers; this is why we have limited our discussion here to isolated points of the spectrum).

(2) If $v\in X_1$, $G_v$ has a simple pole at $0$ and for other $v$'s either a zero or a non-zero finite value at $0$.

(3) A direct analysis of the possibilities proves that if $e=(vw)$ with both points in $X_1$, then $m_e$ has a finite non-zero value at $0$ so, by \eqref{2.5}, $Q_e$ has a simple pole.

(4) A direct analysis of the possibilities proves that if $e=(vw)$ with $v\in X_1, w\in\partial X_1$, then $m_e(0)=0$ and $m_{\check{e}}$ has a pole at $0$ so $Q_e$ has a finite, non-zero value at $0$ (since $m_em_{\check{e}}$ has a negative value at $0$ so the denominator in the first equality in \eqref{2.5} is non-varnishing) and $G_w(0)=0$.

(5) A direct analysis of the possibilities proves that if $e=(vw)$ with both points not in $X_1$, then $Q_e$ does not have a pole at $0$ so by (3) and (4), $Q_e$ has a pole at zero if and only if both endpoints lie in $X_1$.

(6) The final equalities in \eqref{2.5} show that if $e=(vw)$ and $Q_e(0)=0$, then neither $v$ nor $w$ can lie in $X_1$.  It also shows that if $m_e$ has a pole at $0$, and neither $v$ nor $w$ lies in $X_1$, then $Q_e(0)=0$.  It follows that for such $e$'s, $Q_e$ has a double $0$ at $0$ if both $m_e$ and $m_{\check{e}}$  have poles there (by the first equality in \eqref{2.5}) and has a simple pole at $0$ if exactly one of them has a pole.  Thus for such $e$'s, we can count poles of $m_e$ rather than zeros of $Q_e$ so long as we run $e$ through all of $\tilde{E}$.

(7) It follows from \eqref{2.2}, that if $G_u(0)=0$, then at least one $m_f$ with $u=\sigma(f)$ has a pole, and because poles all have negative residues, the converse is true.  A careful analysis shows that if $m_f$ with $\sigma(f)=u$ has a pole and if $e\ne f$ also has $\sigma(e)=u$ and with $\tau(e)\notin X_1$, then $m_e$ is not infinite.  This means if that there is a $1-1$ correspondence between $v\in X_0$ with $G_v(0)=0$ and those $e$ with $m_e$ having a pole with $\sigma(e)\in X_0$.  By the argument in (5), it also says that if $\sigma(e)\in\partial X_1, \tau(e)\notin X_1$, then $m_e$ does not have a pole.  These two conclusions show that the number of zeros of $G_u(0)$ with $u\in X_0$ exactly cancels the number of zeros of $Q_e$, for those $e$ with no ends in $X_1$.

These observations prove the Aomoto index theorem.

We remark that the earlier proofs of Aomoto's theorem (\cite{AomotoPoint, BGM}) show that $X_1$ is a forest (disjoint union of trees) which allows one to prove that the index is also equal to $ccX_1(\lambda)-\#(\partial X_1(\lambda))$ where $ccX_1(\lambda)$ is the number of connected components of $X_1(\lambda)$.  So long as we use the formula \eqref{1.4}, we don't need to prove the forest result.

\section{Anderson Model on a Tree} \lb{s5}

In this final section we will note that the ideas of Section 2 also imply results for the Anderson model on a tree, a subject with considerable work in both the physics \cite{ATA, Derr, Mill} and mathematical physics \cite{Ai1, Ai2, Ai3, Froese, Klein} literatures.  One fixes a strictly positive integer, $d$, and considers a Jacobi matrix on the homogeneous tree of degree $d$. The $a$'s and $b$'s are both given by independent identically distributed (separately for $a$ and $b$) random variables (for technical simplicity, we suppose the supports of the distributions are bounded).  Most commonly the distributions of the $a$'s set them to be identically one but that doesn't affect anything in our arguments.

For us, as for Klein~\cite{Klein}, the density of states is given by the expectation of the spectral measure over the ensemble of random Hamiltonians.  By taking expectations of  \eqref{1.6} and \eqref{2.10}, we prove that
\begin{equation}\label{5.1}
  \int \log(t-z)\, dk(t) = \left(\tfrac{d}{2}-1\right)\bbE(\log(G_u)) - \tfrac{d}{2}\bbE(\log(m_e))
\end{equation}
In case $d=2$ this is what follows from the Thouless formula and \cite[(1.7)]{SiKotani} so this is sort of a half-Thouless formula.  We are currently studying possible applications of \eqref{3.1}.

\ifnum \value{pnas}>1 

\section{PNAS Specific Information} \lb{s6}

\textbf{ORCID}:\\
J. Breuer: 0000-0002-2765-910X \\
J. Garza Vargas: 0000-0001-6258-0600\\
B. Simon: 0000-0003-2561-8539

\acknow{Research of JBr, JGV, ES and BS supported in part by Israeli BSF Grant No. 2020027. Research of JBr and ES supported in part by Israel Science Foundation Grant No. 1378/20. Research of JGZ supported in part by NSF FRG Award 1952777 and Caltech Carver Mead New Adventures Fund under the aegis of Joel Tropp's award, and Caltech CMI Postdoctoral Fellowship}

\showacknow 

\bibsplit[18]

\fi



\end{document}